\documentclass[12pt]{article}
\usepackage{amsfonts}
\usepackage{amsmath}
\usepackage{amssymb}
\usepackage{amsthm}
\usepackage{graphicx}
\usepackage{latexsym}
\makeatletter
\renewcommand{\@seccntformat}[1]
{\csname the#1\endcsname.\enspace}
\makeatother
\newtheorem{theorem}{Theorem}
\newtheorem{lemma}{Lemma}
\newtheorem{remark}{Remark}
\newtheorem{corollary}{Corollary}
\newtheorem{example}{Example}

\newtheorem{assumption}{Assumption}

\begin{document}
\large
\begin{center}
{\bf  On shrinkage estimation for balanced loss functions \footnote{\today}}   
\end{center}
\normalsize
%begin{tabular}
\begin{center}
{\sc \'Eric Marchand$^{a}$ \& William E. Strawderman$^{b}$} \\
    
{\it a  Universit\'e de
    Sherbrooke, D\'epartement de math\'ematiques, Sherbrooke Qc,
    CANADA, J1K 2R1 \quad (e-mail: eric.marchand@usherbrooke.ca) } \\

{\it b  Rutgers University, Department of Statistics, 501 Hill
Center, Busch Campus, Piscataway, N.J., USA, 08855 \quad (e-mail:
straw@stat.rutgers.edu) }
\end{center}
\vspace*{0.2cm}
\begin{center}
{\sc Summary} \\
\end{center}
\vspace*{0.2cm}
The estimation of a multivariate mean $\theta$ is considered under natural modifications of balanced loss function of the form: (i) $\omega \, \rho(\|\delta-\delta_0\|^2) + (1-\omega) \, \rho(\|\delta-\theta\|^2) $,   and (ii)  $\ell \left( \omega \, \|\delta-\delta_0\|^2 + (1-\omega) \, \|\delta-\theta\|^2 \right)\,$, where $\delta_0$ is a target estimator of $\gamma(\theta)$.  After briefly reviewing known results for original balanced loss with identity $\rho$ or $\ell$, we provide, for increasing and concave $\rho$ and $\ell$ which also satisfy a completely monotone property, Baranchik-type estimators of $\theta$ which dominate the benchmark $\delta_0(X)=X$ for $X$ either distributed as multivariate normal or as a scale mixture of normals.  Implications are given with respect to model robustness and simultaneous dominance with respect to either $\rho$ or $\ell$.

\vspace*{0.5cm}

 
\normalsize 
\section{Introduction}

Balanced loss functions and their role in estimation have captured the interest of many researchers over the years since Arnold Zellner (Zellner, 1994) proposed their use in a regression framework.    Balanced loss functions are appealing as they combine proximity of a given estimator $\delta$ to both a target estimator $\delta_0$ and the unknown parameter $\theta$ which is being estimated.   They relate conceptually to methods for combining estimators (e.g., Judge \& Mittlehammer, 2004), as well as penalized least-squares estimation.   The study of balanced loss functions has frequently been cast in a regression framework (e.g., Hu \& Peng, 2011, and the references therein), but it also has arisen or  related to credibility theory, finance, sequential estimation, etc (Baran \& Stepie\'{n}-Baran, 2013; Zhang \& Chen, 2018).  In Zellner's framework,  the target estimator was least-squares, but such a target can be viewed more broadly (e.g., Jafari Jozani et al., 2006, 2014). 

To a large extent, findings in the literature relate to balanced squared error loss
\begin{equation}
\label{blfsquared}
L_{\omega}(\theta, \delta) \, = \, \omega \, \|\delta-\delta_0\|^2 \, + \, (1-\omega) \, \|\delta-\gamma(\theta)\|^2 \,,
\end{equation}
where for an observable $X \sim f_{\theta}$, $\gamma(\theta) \in \Gamma \subset \mathbb{R}^d$, $\delta_0(X)$ is a target estimator of $\gamma(\theta)$,  $\omega \in [0,1)$ is the weight given to the proximity of $\delta$ to $\delta_0$, and  $\delta(X)$ is a given estimator of $\gamma(\theta)$.  In such cases, as presented by Jafari Jozani et al. (2006), as well as Dey et al. (1999),  Bayesian estimation as well as the frequentist risk performance under balanced loss $L_{\omega}$ with $\omega>0$ relate precisely to corresponding features under unbalanced loss (i.e., squared error loss) $L_0$ (see Theorem \ref{lemmajjmp2006}).  
For instance, given a prior $\pi$ and a corresponding Bayes estimator $\delta_{\pi,0}(X)$ under loss $L_0$, the corresponding Bayes estimator under balanced loss $L_{\omega}$ is simply given by $ (1-\omega) \, \delta_{\pi}(X) + \omega \delta_{0}(X)$.  Such relationships are reviewed and briefly illustrated in Section 2.

In contrast, much less is known for the following two natural alternatives or modifications to loss (\ref{blfsquared}):

\begin{equation}
\label{lossrho}
\omega \, \rho(\|\delta-\delta_0\|^2) + (1-\omega) \, \rho(\|\delta-\gamma(\theta)\|^2)
\end{equation} 
and
\begin{equation}
\label{lossell}
 \ell \left( \omega \, \|\delta-\delta_0\|^2 + (1-\omega) \, \|\delta-\gamma(\theta)\|^2 \right)\,,
\end{equation}
with $0 \leq \omega < 1$, $\rho(\cdot) \geq 0$, and $\ell(\cdot) \geq 0$.
Balanced loss functions of the type (\ref{lossrho}) were considered by Jafari et al. (2012).   They provided Bayesian estimators as well as other type of posterior risk analysis.  However, for both losses (\ref{lossrho}) and (\ref{lossell}), there seems to be no significant known finding for frequentist risk analysis, such as the earlier results for balanced squared error loss.   

The objective of this paper is to try to fill such gaps.  To achieve this, we focus on the multivariate normal case $X \sim N_d(\theta, \sigma^2 I_d)$, as well as scale mixture of normals as defined in (\ref{representation-mixture}), the target estimator $\delta_0(X)=X$; and the objective of improving on $X$.  The latter is the maximum likelihood estimator and also is minimax for losses (\ref{lossrho}) as elaborated upon at the outset of Section 3C.  We obtain various sufficient conditions for dominance for both losses (\ref{lossrho}) and (\ref{lossell}).   These apply for interesting subclasses of concave $\rho$'s and $\ell$'s respectively, which are also completely monotone.   Shrinkage estimation for multivariate normal models, and more generally spherically symmetric and elliptically symmetric models, has had a long, rich and influential history (e.g., Fourdrinier et al., 2018).   The use of a concave loss as well as concave versions of (\ref{lossrho}) and (\ref{lossell}), is quite appealing, and has motivated previous shrinkage estimation work such as Brandwein \& Strawderman (1980, 1991), Brandwein et al. (1993), and Kubokawa et al. (2015),  among others.

The paper is organized as follows.  We collect some preliminary definitions and results in Section 2.1, before reviewing and illustrating frequentist risk and 
Bayesian analysis results in Section 2.2 applicable to balanced squared-error loss $L_{\omega}$.  In Sections 3 and 4, we provide conditions for a Baranchik-type estimator to dominate $\delta_0(X)=X$ under loss functions (\ref{lossrho}) and (\ref{lossell}) respectively (i.e., Theorems \ref{dominancesection3} and \ref{theoremdominanceforl}).  In both cases, the proofs are unified with respect to choice of model and loss, the former with respect to the underlying normal mixture and the latter with respect to the choice of $\rho$ or $\ell$ for the balanced loss.   Implications are given in terms of robustness and simultaneous dominance (i.e., Corollary  \ref{corollart-simultaneous-dominance}).   
Finally, we make use of various techniques and properties relative to concave functions, completely monotone functions, superharmonic functions, and spherically symmetric distributions.

\section{Preliminary results and the balanced squared-error loss case}

\subsection{Preliminary definitions and properties}

We assemble here some definitions and properties useful throughout the manuscript.  The estimators studied below are based on spherically symmetric distributions $X \sim f(\|x-\theta\|^2), x, \theta \in \mathbb{R}^d$, which are scale mixtures of normals.  Such distributions admit the representation
\begin{equation}
\label{representation-mixture}
X|V \, \sim \, N_d(\theta, V I_d), \; V \sim g\,, 
\end{equation}
and include many familiar examples such as Normal, Student, Logistic, Laplace, Exponential power (with $f(t)=t^s, 0<s<1$), among others.   Other than moment finiteness conditions and the restriction to $d \geq 3$ or $d \geq 4$ dimensions, the applicability of our dominance findings will not require any further specific assumptions on $f$.

A key characterization and property, which brings into play completely monotone functions, is given by the following result (see, e.g., Feller, 1966; Berger, 1975; etc.).

\begin{lemma} \label{lemma-completelymonotone}
\begin{enumerate}
\item[ {\bf (a)}] A density of the form $f(\|x-\theta\|^2), x,\theta \in \mathbb{R}^d$ is a scale mixture of normals if and only if $f(\cdot)$ is completely monotone, i.e.,
$(-1)^{n} f^{n}(t) \geq 0$, for $ n = 1,2, \ldots$, and $t \in  \mathbb{R}_+ \,$.

\item[ {\bf (b)}] The product of two completely monotone functions is completely monotone.

\end{enumerate}
\end{lemma}

The dominance findings of Sections 3 and 4 relate to Baranchik-type estimators of $\theta$ defined and denoted throughout as:

\begin{equation}
\label{baranchik}
\delta_{a, r(\cdot)}(X) \, = \, \left(1 \, - \, \frac{a r(\|X\|^2)}{\|X\|^2}  \right) \, X\,,
\end{equation}
with $a>0$, and the conditions

\begin{equation}
\label{conditions}
 0 \leq r(\cdot) \leq 1\,, \, r(\cdot) \neq 0\,,\, r'(\cdot) \geq 0\,,\, \hbox{ and } (d/dt) (r(t)/t) \leq 0\, .
\end{equation}
  These include James-Stein estimators with constant $r(\cdot)$, and $r(t)=(d-2) \, \sigma^2$ in the original $X \sim N_d(\theta, \sigma^2 I_d)$ case.

\subsection{Balanced squared-error loss}

%As mentioned above, there exists a strong link between Bayesian posterior analysis and frequentist risk analysis between balanced loss $L_{\omega}$ (\ref{blfsquared}) with $\omega>0$ and squared error loss $L_0$.

We review here, for Bayesian inference and frequentist risk analysis, relationships between balanced loss $L_{\omega}$ and its unbalanced counterpart $L_0$.  Such results appear in Dey et al. (1999), as well as in Jafari Jozani et al. (2006).  For the former, the findings apply to a multivariate normal model $X \sim N_d(\theta, \sigma^2 I_d)$ and $\delta_0(X)=X$, while the latter work relates to a more general model $X \sim f_{\theta}$ and target estimator $\delta_0$. Some of the results will serve in later sections, but they are exposed here also to illustrate the facility in which Bayesian analysis and frequentist risk evaluations for $L_{\omega}$ follow from corresponding results for squared-error loss $L_0$.

The following Lemma \ref{lemmajjmp2006} will be used in Section 4 for the analysis of losses in (\ref{lossell}), but is presented here as it serves to link the frequentist risk  under loss $L_{\omega}$ to the risk under squared error loss $L_{0}$, as presented in Corollary \ref{corollarysel}.    
To facilitate the presentation that follows, we denote the difference in losses $L_{\omega}$ between estimates $\delta_0(x) + (1-\omega) \, g(x)$ and $\delta_0(x)$ as

\begin{equation}
\label{Delta}
\Delta_{\omega}(\theta, g) \, = \, L_{\omega} \left(\theta, \delta_0 + (1-\omega) g\right) \; - \; L_{\omega}\left(\theta, \delta_0\, \right).
\end{equation}

\begin{lemma}
\label{lemmajjmp2006}
Let $X \sim f_{\theta}$.  For the problem of estimating $\gamma(\theta)$ under balanced loss $L_{\omega}$ (as in (\ref{blfsquared})),  we have
$\Delta_{\omega}(\theta, g) \, = \, (1-\omega)^2 \; \Delta_{0}(\theta, g) \, $. 
\end{lemma}
{\bf Proof.}
A decomposition of (2.4) yields
\begin{eqnarray*}
\Delta_{\omega}(\theta, g) \, & = & \, \omega \, \| \delta_0 + (1-\omega) g - \delta_0\|^2 \; + \; (1-\omega) \, \|\delta_0 + (1-\omega) g - \gamma(\theta) \|^2 \, - \,
(1-\omega) \|\delta_0-\gamma(\theta) \|^2 \\
\, & = & \, (1-\omega)^2 \, \|g\|^2 \, + \, 2 (1-\omega)^2 g^\top (\delta_0-\gamma(\theta)) \\
\, & = & \, (1-\omega)^2 \, \left( \|\delta_0 + g - \gamma(\theta)\|^2  
- \, \|\delta_0 - \gamma(\theta)\|^2 \right)  \\
\, & = & \,  (1-\omega)^2 \; \Delta_{0}(\theta, g). \;\;\;\;\;\;\;\;\;\;
\;\;\;\;\;  \qed
\end{eqnarray*}

In terms of the frequentist risk $R_{\omega}$ associated with loss $L_{\omega}$,
given by $R_{\omega}(\theta, \delta) \, = \, \mathbb{E} \, \{L_{\omega} \left(\theta, \delta(X) \right) \,\}$ for an estimator $\delta(X)$ of $\gamma(\theta)$,
the following general result follows from Lemma \ref{lemmajjmp2006}.
   
\begin{corollary} 
%(Jafari Jozani et al., 2006)
\label{corollarysel}
Let $X \sim f_{\theta}$ and consider the problem of estimating $\gamma(\theta)$.  The estimator $\delta_{1,\omega}(X)\,=\,\delta_0(X) + (1-\omega) g_1(X)$ dominates 
$\delta_{2,\omega}(X)\,=\,\delta_0(X) + (1-\omega) g_2(X)$ under loss $L_{\omega}$ if and only if $\delta_{1,0}(X) \, = \, \delta_0(X) + g_1(X)$ dominates $\delta_{2,0}(X) \, = \, \delta_0(X) + g_2(X)$ under squared error loss $L_{0}$. 
\end{corollary}
{\bf Proof.}   We have
\begin{eqnarray*}
R_{\omega}(\theta, \delta_{1,\omega}) \, - \, 
R_{\omega}(\theta, \delta_{2,\omega}) & = & \mathbb{E}_{\theta} \{\Delta_{\omega}(\theta, g_1(X)) - \Delta_{\omega}(\theta, g_2(X))  \} \\
\, & = & \, (1-\omega)^2 \,\, \mathbb{E}_{\theta} \{\Delta_{0}(\theta, g_1(X)) - \Delta_{0}(\theta, g_2(X))  \} \; \;\;  (\hbox{Lemma \ref{lemmajjmp2006})} \\
\, & = & \, (1-\omega)^2 \,\, \{ R_{0}(\theta, \delta_{1,0}) \, - \, 
R_{0}(\theta, \delta_{2,0}) \},
\end{eqnarray*}
which establishes the result.  \qed

Now, turning to Bayesian inference, we have an equally simple relationship between balanced loss $L_{\omega}$ and its unbalanced counterpart $L_0$.   
More precisely, the following well-known result conveniently expresses the Bayes estimator $\delta_{\pi, \omega}$ under $L_{\omega}$ for $\omega >0$ in terms of the Bayes estimator $\delta_{\pi, 0}$ under $L_0$, given of course by $\delta_{\pi, 0}(X) \, = \, \mathbb{E} (\gamma(\theta)|X)$.

\begin{theorem}
\label{theorem-bayes-sel}
For $X \sim f_{\theta}$ and a prior $\theta \sim \pi$ for which $Cov(\gamma(\theta)|x) $ exists for all $x$, the Bayes estimator $\delta_{\pi, \omega}$ of 
$\gamma(\theta)$ under loss $L_{\omega}$ is given by $\delta_{\pi, \omega}(X) \, = \, \omega \, \delta_0(X) \,+ \, (1-\omega) \, \delta_{\pi,0}(X)\,. $
\end{theorem}
{\bf Proof.}  Write $\delta_{\pi, \omega}(x) \, = \, \delta_0(x) \, + \, (1-\omega) \, g_{\pi,\omega}(x) $ for $0 \leq \omega < 1$.  By definition of the Bayes estimate, we thus have
\begin{eqnarray*}
g_{\pi,\omega}(x) \, & = & \,  \underset{g}{\arg\min} \, \mathbb{E} \{L_{\omega}\left(\theta, \delta_0(x) \, + \, (1-\omega) g \right) \, | \, x  \}  \\
\, & = & \,  \underset{g}{\arg\min}  \, \mathbb{E} \{\Delta_{\omega}(\theta,  g ) \, | \, x  \} \\
\, & = & \,  \underset{g}{\arg\min}  \, \left[(1-\omega)^2 \, \mathbb{E} \{\Delta_{0}(\theta,  g ) \, | \, x  \} \right] \; \;\;  (\hbox{Lemma \ref{lemmajjmp2006})} \\  \, & = & \,  \underset{g}{\arg\min} \, \mathbb{E} \{\|\delta_0(x) \, + \, g \, - \, \gamma(\theta) \|^2 \, x  \} \\  \, & = & \, \mathbb{E} \{\gamma(\theta)|x \} \, - \, \delta_0(x)\,.  
\end{eqnarray*}
From this, the result follows as $\delta_{\pi, \omega}(x) \, = \, \delta_0(x) \, + \, (1-\omega)  \left[\mathbb{E} \{\gamma(\theta)|x \} \, - \, \delta_0(x)\,\right] \, = \, \omega \, \delta_0(x) \,+ \, (1-\omega) \, \delta_{\pi,0}(x)\,.$
\qed

Now, combining the last two results leads to the following Bayes dominance result.

\begin{corollary}
\label{corollarybayesdominance}
For $X \sim f_{\theta}$, a prior $\theta \sim \pi$, and the problem of estimating $\gamma(\theta)$, the Bayes estimator $\delta_{\pi, \omega}(X)$ dominates an estimator $\delta_0(X)+ (1-\omega) g_2(X)$ under $L_{\omega}$ if and only if the Bayes estimator $\delta_{\pi, 0}(X)\, = \, \mathbb{E} (\gamma(\theta)|X) $ dominates $\delta_0(X)+ g_2(X)$ under squared-error loss $L_0$.
\end{corollary}

Several examples can be found in the literature, namely among the references mentioned above.    We do provide at the end of this subsection Example \ref{examplesel} as an illustration. Before doing so, we briefly address the issue of minimaxity, where relationships between balanced and unbalanced losses are not as immediate (e.g., Jafari Jozani et al., 2006).  One situation though does simplify, namely the case where the target estimator $\delta_0(X)$ is itself minimax under the unbalanced loss.  Moreover, the following result (Jafari Jozani et al., 2012, Theorem 4) holds in general for losses (\ref{lossrho}).   As discussed at the outset of Section 3C, this will serve to guarantee that the dominating estimators of Theorem \ref{dominancesection3} are themselves minimax.     

\begin{theorem}
\label{theorem-minimax}  Let $X \sim f_{\theta}$ and consider the problem of estimating 
$\gamma(\theta)$ under loss (\ref{lossrho}).  Suppose that the estimator $\delta_0(X)$ is minimax under unbalanced loss $\rho(\|\delta-\gamma(\theta)\|^2)$.  Then, $\delta_0(X)$ is also minimax under loss (\ref{lossrho}) for all $0<\omega<1$.    
\end{theorem}
{\bf Proof.}
Let $R_w$ denote the frequentist risk under loss (\ref{lossrho}).  Since $R_{\omega}(\theta, \delta_0) \, = \, (1-\omega) R_0(\theta, \delta_0)$, the result is immediate.  \qed

\begin{example}
\label{examplesel}
We consider the classical problem of estimating a multivariate normal mean and 
illustrate how known Stein estimation results applicable (e.g., Stein 1981; Strawderman, 2003)  to squared-error loss $L_0$ translate to balanced loss $L_{\omega}$.  Let $X \sim N_d(\mu, \sigma^2 I_d)$ and $S^2 \sim \sigma^2 \chi^2_k$ be independently distributed with $d \geq 3, k \geq 1$. Set $\theta=(\mu, \sigma^2)$ and consider estimating $\tau(\theta)=\mu$ under balanced loss $L_{\omega}$ with target estimator $\delta_0(X)=X$, which is minimax under $\frac{L_0(\theta, \delta)}{\sigma^2}$, and thus minimax under loss $\frac{L_{\omega}(\theta, \delta)}{\sigma^2} \, = \, \frac{ \omega \, \|\delta-X\|^2 + (1-\omega) \, \|\delta-\theta\|^2}{\sigma^2}$  by virtue of Theorem \ref{theorem-minimax}.

\begin{enumerate}

\item[ {\bf (A)}]  For known $\sigma^2$, any estimator of the form $\delta(X)\, = \, X \, + \, \sigma^2 \, g(X)$ such that $g$ is weakly differentiable, $\mathbb{E}_{\theta} \|g(X) \|^2 < \infty$, and $\|g(X)\|^2 + 2 \, div (g(X)) \leq 0$ a.e., dominates $\delta_0(X)$ under loss $L_0$.  It thus follows from Corollary \ref{corollarysel} that $X + (1-\omega) \, \sigma^2 \, g(X)$ dominates
$\delta_0(X)$ under loss $L_{\omega}$ for such $g$'s.  Such dominating estimators include the James-Stein estimator with $g(t)=-(d-2) t/\|t\|^2$, as well as Baranchik type estimators $\delta_{a, r(\cdot)}(X)$ in (\ref{baranchik}), with $0 < a < 2(d-2) (1-\omega) \sigma^2$ and conditions (\ref{conditions}) on $r(\cdot)$.

\item[ {\bf (B)}]   For unknown $\sigma^2$, with the same conditions on $g$, estimators of the form $X \, + \, \frac{S^2}{k+2} \, g(X)$ dominate $X$ under loss $L_0$.  Again, it follows immediately from Corollary \ref{corollarysel} that estimators $X \, + \, (1-\omega) \, \frac{S^2}{k+2} \, g(X)$ dominate $X$ under balanced loss $L_{\omega}$.

\item[ {\bf (C)}]   For known $\sigma^2$,  Bayes estimators $\delta_{\pi, \omega}(X)$ under balanced loss $L_{\omega}$ and associated with prior density 
$\pi(\theta)$, Theorem \ref{theorem-bayes-sel} along with a well-known representation for $\delta_{\pi, 0}(X)$ tell us that
\begin{eqnarray*}
\delta_{\pi, \omega}(X) \, &=& \, \omega X \, + \, (1-\omega) \delta_{\pi, 0}(X)
 \\
 \, &=& \, \omega X \, + \, (1-\omega) (X + \sigma^2 \frac{\nabla m(X)}{m(X)}) \, = \,  X \, + \, (1-\omega) \, \sigma^2 \frac{\nabla m(X)}{m(X)}\,,
\end{eqnarray*}
where $m(X) \, = \, (2\pi \sigma^2)^{-1} \, \int_{\mathbb{R}^d} e^{- \frac{1}{2 \sigma^2} \|X-\theta\|^2} \pi(\theta) \, d\theta$ is the marginal distribution of $X$.   By virtue of Corollary \ref{corollarybayesdominance}, the estimator 
$\delta_{\pi, \omega}(X)$ dominates $X$ under loss $L_{\omega}$ if and only if
$\delta_{\pi, 0}(X)$ dominates $X$ under loss $L_{0}$.  With the superharmonicity of either $\pi(\cdot), m(\cdot)$ or $\sqrt{m(\cdot)}$ a sufficient condition for $\delta_{\pi, 0}(X)$ to dominate $X$ under loss $L_{0}$ (e,g. Strawderman, 2003), we thus infer that either of these conditions imply that $\delta_{\pi, \omega}(X)$ dominates $X$ under balanced loss $L_{\omega}$. 
\end{enumerate}
\end{example}

\section{Risk analysis for loss $\omega \rho(\|\delta-X\|^2) + (1-\omega) \rho(\|\delta-\theta\|^2)$}

%\subsection*{A.  The model}

\subsection*{A.  The loss function}

For a model (\ref{representation-mixture}), we evaluate the frequentist risk performance of an estimator $\delta(X)$ of $\theta$ under the balanced loss

\begin{equation}
\label{loss-section3}
L_{\omega,\rho}(\theta, \delta) \, = \, \omega  \,  \rho(\|\delta-X\|^2) + (1-\omega)  \,  \rho(\|\delta-\theta\|^2)\,, 0 \leq \omega <1 \,, 
\end{equation}
which incorporates the target estimator $\delta_0(X)=X$.  For the function $\rho$, we assume the following throughout this section:

\begin{assumption}
\label{assumptionrho}
$\rho(0)=0 \,,\, 0 <\rho'(0) < \infty  $, and $\rho'$ is completely monotone on $\mathbb{R}_+$, i.e., $(-1)^{k} \, \rho^{(k+1)}(t) \geq 0$ for $t>0$ and for $k=0,1, \ldots$. 
\end{assumption}

Examples of loss functions $L_{\omega,\rho}$ for which $\rho$ satisfies Assumption \ref{assumptionrho}, other than $\rho(t)=t$, include: {\bf (i)} $\rho(t) = 1 - e^{-t/\alpha}$ with $\alpha>0$, {\bf (ii)} $\rho(t) \, = \, \log(1+t)$, {\bf (iii)} $\rho(t) = (1+ t/\gamma)^{\beta}$ with $\gamma >0, \beta \in (0,1)$, and  {\bf (iv)}  cases $\rho(t) \, = \, z(0) - z(t)$ with $z$ being completely monotone such as $\rho(t) \, = \, r^2t/(rt+1)$ with $r>0$.  Case {\bf (i)} is known as reflected normal loss, while examples {\bf (iv)} represent a broader class of bounded losses.   $L^{\beta}$ losses with $\rho(t)=t^{\beta}$, $\beta \in (0,1)$, represent concave choices, but such $\rho$'s do not satisfy the finiteness assumption on $\rho'(0)$.

\subsection*{B.  Further technical results}

We now expand on various technical results which are pivotal to the risk analysis in Subsection {\bf 3C}.

\begin{lemma}
\label{lemmafvs.f*}

Consider $X \sim f(\|x-\theta\|^2), x, \theta \in \mathbb{R}^d$, admitting representation (\ref{representation-mixture}) with mixing variable $V$, and $\rho$ satisfying Assumption \ref{assumptionrho}.  Let $Y \sim f^*(\|y-\theta\|^2) \, = \, \rho'(\|y-\theta\|^2) \, f(\|y-\theta\|^2)/K  \,; y \in \mathbb{R}^d$; with $K=\mathbb{E}_0 (\rho'(\|X\|^2))$.  Then,

\begin{enumerate}
\item[{\bf (a)}]  The distribution of $Y$ admits a scale mixture of normals representation  
\begin{equation}
\label{distributionY}
Y|W \sim N_d(\theta, W I_d)\,, \, \hbox{ with } W \sim h\,;
\end{equation}
\item[{\bf (b)}]   Moreover, the distribution of $W$ is stochastically smaller than the distribution of $V$.
\end{enumerate}

\end{lemma}
{\bf Proof.}  First observe that $f^*$ is a density since $K \, = \, \mathbb{E}_0 (\rho'(\|X\|^2)) \leq \rho'(0) < \infty$.   Part {\bf (a)} thus follows from Lemma \ref{lemma-completelymonotone}.  For part {\bf (b)}, given that $\rho'$ and $f$ are completely monotone, they are representable as Laplace transforms (Lemma \ref{lemma-completelymonotone}):
\begin{eqnarray*}
f(t) \, = \, K_1 \, \int_0^{\infty}  e^{-t/2v} \, dG(v)\,,\\
\rho'(t) \, = \, K_2 \, \int_0^{\infty}  e^{-t/2\tau} \, dH(\tau)\,,
\end{eqnarray*}
for $t \in \mathbb{R}^d$.  From this, we have for $t \in \mathbb{R}^d$
\begin{eqnarray*}
f^*(t) \, & = & \frac{K_1 K_2}{K} \int_0^{\infty} \int_0^{\infty} e^{-\frac{t}{2} (\frac{1}{v} + \frac{1}{\tau})} \, dG(v) \, dH(\tau) \;,\\
\,  & = & K_3 \, \int_0^{\infty} e^{-t/2w} dM(w) \,.
\end{eqnarray*}
Interpreting in terms of scale mixture of normals, we have for $Y \sim f^*(\|y-\theta\|^2)$ representation (\ref{distributionY}) with 
$W=^d\frac{\tau V}{\tau + V}$.  
Finally, from this, we have $\mathbb{P}(W \leq s) \, \geq \,\,\mathbb{P}(V \leq s),$ for all $s>0$ and the result follows.
\qed

The two lemmas that follow, which we will require, rely partly on properties of  superharmonic functions.   We recall that a continuous function $g: \mathbb{R}^d \to \mathbb{R}$ is superharmonic if and only if: at all $t_0 \in  \mathbb{R}^d$ and $r>0$, the average of $g$ over the surface of the sphere, centered at $t_0$ of radius $r$, $S_r(t_0)=\{t \in \mathbb{R}^d: \|t-t_0\| =r  \}$ is less or equal than $g(t_0)$.  For twice differentiable $g$, the superharmonicity of $g$ is equivalent to its Laplacian being less or equal to $0$, i.e., $\Delta \, g \leq 0$ with $\Delta \, g \, = \, \sum_{i=1}^d \frac{\partial^2}{\partial t_i^2} g(t)$.

\begin{lemma}
\label{lemmatechnical}
Let $Z \sim N_d(0,I_d)$ with $d \geq 3$ and let $T=\|\alpha Z + \theta\|^2$ with $\alpha>0$ and $\theta \in \mathbb{R}^d$.  Then, we have the following:
\begin{enumerate}
\item[{\bf (a)}]
$\epsilon(\alpha) \, = \, \mathbb{E}_{\alpha} (\frac{1}{T})$ is decreasing in $\alpha$ for $\alpha>0$;
\item[{\bf (b)}]  $\mathbb{E} \, g(\alpha Z + \theta)$ is non-increasing in $\alpha$ provided that $g(\cdot) >0$ and $g$ is superharmonic.
\end{enumerate}

\end{lemma}
{\bf Proof.}    The proof of part {\bf (a)} is relegated to an Appendix.  For part {\bf (b)}, first denote $U_m, m>0,$ as a random vector uniformly distributed on the sphere $S_m(0)$ centered at $0$ of radius $m$.  
It suffices to show that $\beta(\alpha,r) = \mathbb{E} \, \left( g(\alpha Z+\theta)|\|Z\|=r \right)$ is for all $r>0$ decreasing in $\alpha$.   Since $(Z|\|Z\|=r) \, \sim \, U_r$ independently of $\|Z\|$, we have $\beta(\alpha,r) \, = \, \mathbb{E} \, \left\lbrace g(U_{\alpha r}+\theta)\right\rbrace $.
Since, for a superharmonic function, the sphere mean is decreasing in the radius (see, e.g., Fourdrinier et al. 2018, Theorem 7.4), we infer that $\beta(\alpha,r)$ is decreasing in $\alpha$, which concludes the proof.    \qed

\begin{lemma}
\label{lemmatechnical2}
Let $\theta \in \mathbb{R}^d$, $a>0$, and $\rho$ satisfy Assumption \ref{assumptionrho}.   Consider $X \sim f(\|x-\theta\|^2)$, $Y \sim f^*(\|y-\theta\|^2)$ as in (\ref{representation-mixture}) and Lemma \ref{lemmafvs.f*}, respectively.    

\begin{enumerate}
\item[{\bf (a)}]  For $d \geq 3$, we have 
\begin{equation}
\nonumber
\mathbb{E}_{\theta} \left(\rho(\frac{a^2}{X'X}) \right) \, \leq \, \rho'(0) \, \mathbb{E}_{\theta} (\frac{a^2}{X'X}) \, \leq \,  
\rho'(0) \, \mathbb{E}_{\theta} (\frac{a^2}{Y'Y})\,;
\end{equation}

\item[{\bf (b)}]   For $d \geq 4$ and $r: \mathbb{R}^d \to [0,1]$  a twice-differentiable function that is non-decreasing and concave, we have

\begin{equation}
\nonumber
\mathbb{E}_{\theta} \left(\rho(\frac{a^2 \, r^2(\|X\|^2)}{\|X\|^2}) \right) \, 
\leq \, \rho'(0) \, \mathbb{E}_{\theta} \left(\frac{a^2 \, r(\|Y\|^2)}{\|Y\|^2} \right) \,.
\end{equation}
\end{enumerate}
\end{lemma}
{\bf Proof.}  {\bf (a)} The first inequality follows from the inequality
\begin{equation}
\label{rho(t)lessthan}
\rho(t) \, \leq \rho(0) + \rho'(0) t \, = \, \rho'(0) \, t\,,
\end{equation}
which holds since $\rho$ is concave with $\rho(0)=0$.  The second inequality
follows from Lemma \ref{lemmafvs.f*} and part (a) of Lemma \ref{lemmatechnical}.  Indeed,  since $X|V \sim N_d(\theta, V I_d)$, $Y|W \sim N_d(\theta, W I_d)$, we have, with the notation of Lemma \ref{lemmatechnical},  
$\mathbb{E}_{\theta}(\frac{1}{\|X\|^2}) \, = \mathbb{E} (\epsilon (\sqrt{V}))$ and 
$\mathbb{E}_{\theta}(\frac{1}{\|Y\|^2}) \, = \mathbb{E} (\epsilon (\sqrt{W}))$, and the result follows since $\epsilon(\cdot)$ is decreasing  and $W$ is stochastically smaller than $V$.  \\

{\bf (b)}   Defining $Z \sim N_d(0,I_d)$ and denoting $g_0(t) \, = \, \frac{r(\|t\|^2)}{\|t\|^2}\,,\, t \in \mathbb{R}^d$, we have
\begin{eqnarray*}
\mathbb{E}_{\theta} \left(\rho(\frac{a^2 \, r^2(\|X\|^2)}{\|X\|^2}) \right)  
& \leq &  \mathbb{E}_{\theta} \left(\rho(\frac{a^2 \, r(\|Y\|^2)}{\|Y\|^2}) \right)   \\
& \leq &  a^2 \, \rho'(0) \,  \mathbb{E}_{\theta} \left(\frac{ r(\|X\|^2)}{\|X\|^2} \right)  \\
& = &  a^2 \, \rho'(0) \, \mathbb{E} \left(g_0(\sqrt{V} \, Z \, + \, \theta)  \right)  \\
& \leq & a^2 \, \rho'(0) \, \mathbb{E} \left(g_0(\sqrt{W} \, Z \, + \, \theta)  \right)  \\
& = &   a^2 \, \rho'(0) \,  \mathbb{E}_{\theta} \left(\frac{ r(\|Y\|^2)}{\|Y\|^2} \right) \,,
\end{eqnarray*}
where (i) the two equalities follow from the scale mixture representations of $f$ and $f^*$; (ii) the first inequality follows since $\rho$ is non-decreasing and $0 \leq r^2(\|t\|^2) \leq r(\|t\|^2) \leq 1$ for $t \in \mathbb{R}^d$, (iii) the second inequality follows from \eqref{rho(t)lessthan}, and (iv) the third inequality follows from Lemma \ref{lemmafvs.f*}, part {\bf (b)} of Lemma \ref{lemmatechnical}, as in the above proof of part {\bf (a)}, and from the fact that

\begin{equation}
g_0(t) \hbox{ is superharmonic for } d \geq 4,
\end{equation}
provided $r(\cdot)$ is non-negative, non-decreasing, and concave.  Finally, to justify the above, note that, for twice-differentiable $h(\|t\|^2), t \in  \mathbb{R}^d$,
$$  \Delta  \, h(\|t\|^2) \, = \, 2 d h'(\|t\|^2) \, + \, (t't) \, h''(\|t\|^2)\,,$$
so that the choice $h(\|t\|^2)= \, g_0(t) \, = \, \frac{r(\|t\|^2)}{\|t\|^2}$ yields with a little bit of computation
\begin{eqnarray*}
\Delta \left( g_0(t) \right) \, & = & \, 2 \left( (\frac{d-4}{\|t\|^2}) \, \{\|t\|^2 \; r'(\|t\|^2) - r(\|t\|^2) \} \, + 2 r''(\|t\|^2)    \right) \\
\,  & \leq &   0\,
\end{eqnarray*}
since the properties of $r(\cdot)$ imply that $r''(u) \leq 0$ and $r(u) \geq u r'(u)$ for all $u>0$.  \qed

\subsection*{C.  Dominance results}

For balanced loss $L_{\omega, \rho}$ with $\rho$ satisfying Assumption \ref{assumptionrho}, a scale mixture of normals distribution on $X$ with $d \geq 3$, we provide James-Stein and Baranchick-type estimators that dominate $X$.
In such cases, it follows that $X$ is minimax for the unbalanced case $L_{0, \rho}$ with constant risk  $R_0$ (e.g., Kubokawa et al., 2015).  By virtue of Theorem \ref{theorem-minimax}, $X$ is also minimax for balanced loss $L_{\omega, \rho}$.   The following dominance results thus provide dominating estimators which are also minimax under loss $L_{\omega, \rho}$.

\begin{theorem}
\label{dominancesection3}
Consider $X \sim f(\|x-\theta\|^2)$; $x,\theta \in \mathbb{R}^d$; admitting representation (\ref{representation-mixture}), balanced loss function $L_{\omega, \rho}$ as in \eqref{loss-section3} with $\rho$ satisfying Assumption \ref{assumptionrho}. 
\begin{enumerate}
\item[ {\bf (a)}]  If $d \geq 3$, $\delta_a(X) = (1 - \frac{a}{\|X\|^2}) X$ dominates $\delta_0(X)=X$ provided
\begin{equation}
\label{cutoff1}
0< a < \frac{2 (d-2) K (1-\omega) \, \left\lbrace \mathbb{E} (W^{-1}) \right\rbrace^{-1}}{\omega \rho'(0) + (1-\omega) K}\,,
\end{equation}
with $K=\mathbb{E}_0 (\rho'(\|X\|^2))$, and $W$ the mixing variance for $Y \sim f^*(\|y-\theta\|^2)$ as defined in Lemma \ref{lemmafvs.f*}. An equivalent expression for the above dominance condition is
\begin{equation}
\label{cutoff2}
0< a < \frac{2 K^2 (1-\omega) \, \left\lbrace \mathbb{E}_0 \left(\frac{\rho'(\|X\|^2)}{\|X\|^2} \right) \right\rbrace^{-1}}{\omega \rho'(0) + (1-\omega) K}\,.
\end{equation}

\item[ {\bf (b)}]  If $d \geq 4$, a Baranchik-type estimator $\delta_{a,r(\cdot)}(X)$ in (\ref{baranchik})  dominates $\delta_0(X)=X$ provided (\ref{cutoff1}) holds and provided $r(\cdot)$ satisfies conditions (\ref{conditions}).
\end{enumerate}

\end{theorem}
{\bf Proof.}
{\bf (a)}  First, the stated equivalence between \eqref{cutoff1} and \eqref{cutoff2} holds since, on one hand, $$\mathbb{E}_0 (\frac{1}{(\|Y\|^2)} \, = \, \mathbb{E} \{\frac{1}{W} \mathbb{E}_0 (\frac{W}{\|Y\|^2} \,|W) \} = \frac{1}{d-2} \; \mathbb{E} \left(\frac{1}{W}\right)\,, $$
(as $\frac{\|Y\|^2}{W}|W \sim \chi^2_d(0)$ when $\theta=0$)
and, on the other hand,
 $$\mathbb{E}_0 (\frac{1}{(\|Y\|^2)} \, =  \int_{\mathbb{R}^d} \, f^*(\|y\|^2) \, \frac{1}{(\|y\|^2)} \, dy \,  = \, \frac{1}{K} \,  \int_{\mathbb{R}^d} \, \rho'(\|x\|^2) \, f(\|x\|^2) \, \frac{1}{(\|x\|^2)} \, dx \,.$$
 
Second, we have for a difference in risks
\begin{eqnarray*}
\Delta_a(\theta) & = & \, R(\theta, \delta_a) \, - \, R(\theta, \delta_0) \\
\, & = & \, \mathbb{E}_{\theta}  \left[ \omega  \rho(\|\delta_a(X) - X\|^2) \, + \, 
(1-\omega) \, \rho(\| \delta_a(X) - \theta  \|^2) \, - \, (1-\omega) \, \rho(\| X - \theta \|^2)  \right]  \\
\, & = & \, \mathbb{E}_{\theta}  \left[ \omega  \rho(\|\frac{aX}{\|X\|^2}\|^2) \, + \, 
(1-\omega)  \left( \, \rho(\| (1 - \frac{a}{\|X\|^2}) X - \theta  \|^2) \, - \, \, \rho(\| X - \theta  \|^2)  \right) \right] \\
\, & \leq &  \mathbb{E}_{\theta}  \left[ \omega  \rho'(0) \, \frac{a^2}{\|Y\|^2} 
 \, + \, 
(1-\omega) \, K \, \frac{\rho'(\|X-\theta\|^2) }{K} \, \left( \frac{a^2}{\|X\|^2} \, - \, \frac{2 a X' (X-\theta)}{\|X\|^2}  \right) \right] \\
\, & = & \mathbb{E}_{\theta}  \left[ \omega  \rho'(0) \, \frac{a^2}{\|Y\|^2} 
  \, + \,
(1-\omega) \, K \, \left( \frac{a^2}{\|Y\|^2} \, - \, \frac{2 a Y' (Y-\theta)}{\|Y\|^2} \, \right) \right]  \,,
\end{eqnarray*}
where the inequality follows from part {\bf (a)} of Lemma \ref{lemmatechnical2} and the concave function inequality $\rho(t_1) - \rho(t_2) \leq \rho'(t_1) (t_1-t_2)$ for all $t_1,t_2 \geq 0$.   Now, with representation (\ref{distributionY}), by conditioning on $W$, and by the Stein's identity and calculation $\mathbb{E}_{\theta} \left[ (Y-\theta)'\frac{Y}{Y'Y}) |W \right] \, = W \,\mathbb{E}_{\theta} \, \hbox{div} \left(\frac{Y}{Y'Y} \right) \, = \, W \,\mathbb{E}_{\theta} \left(\frac{d-2}{Y'Y}\right)$
(with probability $1$), we obtain

\begin{equation}
\label{Delta_a}
\Delta_a(\theta) \, \leq \,  a \, \mathbb{E}^W \left\lbrace  \mathbb{E}_{\theta}^{Y|W}
\frac{W}{\|Y\|^2} \left( \frac{a \left(\omega \rho'(0) + (1-\omega) K \right)}{W}  - 2 (1-\omega) K (d-2) \right)\,
 \right\rbrace \,.
\end{equation}
By noticing that  $\mathbb{E}_{\theta}^{Y|w}
\left(\frac{W}{\|Y\|^2} \right)$  is increasing in $w>0$, given that 
$\frac{\|Y\|^2}{W} |W \sim \chi^2_d(\lambda=\frac{\|\theta\|^2}{W})$ and $\chi^2_d(\lambda)$ distributions are stochastically increasing in $\lambda$,  we infer from (\ref{Delta_a}) and the covariance inequality (i.e., $\mathbb{E} f_1(W) f_2(W) \leq \mathbb{E} f_1(W) \, \mathbb{E} f_2(W) $ for $f_1(\cdot)$ increasing and 
$f_2(\cdot)$ decreasing) that

\begin{equation}
\label{Delta_a2}
\Delta_a(\theta) \, \leq \, a \, \mathbb{E} \left\lbrace \mathbb{E}_{\theta}^{Y|W}
\frac{W}{\|Y\|^2}
\right\rbrace \, \mathbb{E} \, \left( \frac{a \left(\omega \rho'(0) + (1-\omega) K \right)}{W}  - 2 (1-\omega) K (d-2) \right) \,.
\end{equation} 
From the above, it follows immediately that \eqref{cutoff1} is a sufficient condition for 
$\Delta_a(\theta)$ to be negative for all $\theta$.

{\bf (b)}   The proof is similar to that of part {\bf (a)}.  Using the concave function inequality $\rho(t_1) - \rho(t_2) \leq \rho'(t_1) (t_1-t_2)$, Stein's identity, and part {\bf (a)} of Lemma \ref{lemmatechnical2}, we obtain for the difference in risk
\begin{eqnarray*}
\Delta_a(\theta) & = & \, R(\theta, \delta_{a,r(\cdot)}) \, - \, R(\theta, \delta_0) \\
\, & \leq & a \, \mathbb{E}^W \left\lbrace  \mathbb{E}_{\theta}^{Y|W}
\frac{W \, r(\|Y\|^2)}{\|Y\|^2} \left( \frac{a \left(\omega \rho'(0) + (1-\omega) K \right)}{W}  - 2 (1-\omega) K (d-2) \right)\,
 \right\rbrace  \,. 
\end{eqnarray*}
Now, it is easy to verify that $r(t)/t$ is decreasing in $t>0$ under the given conditions on $r(\cdot)$.  Finally, an application of the covariance inequality leads to an inequality as in \eqref{Delta_a2} with $\mathbb{E}_{\theta}^{Y|W}
\left(\frac{W}{\|Y\|^2} \right)$ replaced by $\mathbb{E}_{\theta}^{Y|W}
\left(\frac{W \, r(\|Y\|^2)}{\|Y\|^2} \right)$.  The result then follows.   \qed

\begin{remark}
From inequality (\ref{Delta_a2}), it also follows that dominance occurs, in both parts {\bf (a)} and {\bf (b)} of Theorem \ref{dominancesection3} for the quantity $a$ equal to the upper cut-off point in (\ref{cutoff1}) (or (\ref{cutoff2})) unless 
$\rho(t)=t$ and $W$ is degenerate, i.e., original balanced loss and the multivariate normal case.
\end{remark}

The proof of Theorem \ref{dominancesection3} is unified with respect to the choice of $\rho$, the coefficient $\omega$ in the balanced loss, and the underlying scale mixture or normals distribution.  To conclude, we point out that the above result can be seen as extensions of Kubokawa et al. (2015), as well as Strawderman (1974), whose results can be seen as particular cases of $\omega=0$ in the former case, and $\omega=0, \rho(t)=t$ in the latter case.

\section{Risk analysis for loss
$ \ell \left( \omega \|\delta-\delta_0\|^2 + (1-\omega) \|\delta-\theta\|^2 \right)$}

The main dominance finding of this section (Theorem \ref{theoremdominanceforl}) relates to a multivariate normal $X \sim N_d(\theta, \sigma^2 I_d)$, and more generally to $X \sim f(\|x-\theta\|^2)$ distributed as a scale mixture of normals as in (\ref{representation-mixture}).
We assess the frequentist risk performance of an estimator $\delta(X)$ of $\gamma(\theta)$ under the balanced loss
\begin{equation}
\label{loss-section4}
L_{\omega, \ell}(\theta, \delta) \, = \,  \ell \left( \omega \, \|\delta-\delta_0\|^2 + (1-\omega) \, \|\delta-\gamma(\theta)\|^2 \right) \,, 0 \leq \omega < 1,
\end{equation}
More specifically, we consider the target estimator $\delta_0(X)=X$ and set  $\gamma(\theta)=\theta$, and our objective is to provide, for $d \geq 3$,  estimators of $\theta$ that dominate $\delta_0(X)=X$ under balanced loss (\ref{loss-section4}) other than Section 2's results for $\ell(t)=t$.
For the function $\ell$, we assume, unless stated otherwise, the following throughout this section:

\begin{assumption}
\label{assumptionell}
$\ell(\cdot) \geq 0$, $\, \ell'(\cdot) > 0 $, and $\ell'$ is completely monotone on $\mathbb{R}_+$, i.e., $(-1)^{k} \, \ell^{(k+1)}(t) \geq 0$ for $t>0$ and for $k=0,1, \ldots$. 
\end{assumption}

Examples of losses $L_{\omega, \ell}$ with $\ell$ satisfying Assumption \ref{assumptionell} include examples {\bf (i), (ii), (iii), (iv)} given for $\rho$ in part {\bf B.} of Section 3, but the cases $\ell(t) = t^{\beta}, \beta \in (0,1),$ are also included here since the assumption $\ell'(0) < \infty$ is not required.

We proceed with a preparatory lemma which exploits the concavity of $\ell$, and which relates the difference in losses $L_{\omega, \ell}$, between estimates 
$\delta_0(x) + (1-\omega) \, g(x)$ and $\delta_0(x)$, to the balanced squared-error loss difference $\Delta_{\omega}(\theta, g)$ in (\ref{Delta}). We therefore define 

\begin{equation}
\label{Deltaell}
\Delta_{\omega,\ell}(\theta, g) \, = \, L_{\omega,\ell} \left(\theta, \delta_0 + (1-\omega) g\right) \; - \; L_{\omega,\ell}\left(\theta, \delta_0\, \right),
\end{equation}
and we have the following.

\begin{lemma}
\label{lemmaDeltaforl}
Let $X \sim f_{\theta}$.  For the problem of estimating $\theta$ under loss (\ref{loss-section4}) with twice-differentiable, increasing, and concave $\ell$, we have
\begin{equation}
\label{Deltaellleq}
\Delta_{\omega,\ell}\,(\theta, g) \, \leq \, (1-\omega)^2 \,\ell' \left\lbrace 
(1-\omega) \, \|\delta_0 - \gamma(\theta) \|^2 \right\rbrace  \, \Delta_0\,(\theta, g)\,.
\end{equation}
\end{lemma}
{\bf Proof.}  The proof uses the fact that $\ell(a+b) \, - \, \ell(a) \leq b \, \ell'(a)$, since $\ell$ is concave, with $a= L_{\omega}(\theta, \delta_0)$ and 
$a+b = L_{\omega}(\theta, \delta_0 + (1-\omega) g)$.  This yields :
\begin{eqnarray*}
\Delta_{\omega,\ell}(\theta, g) \,& = \,& \ell\left(L_{\omega}(\theta, \delta_0+ (1-\omega) g)  \right) \, - \, 
\ell\left(L_{\omega}(\theta, \delta_0)  \right) \\
\, & \leq \, &\ell'\left( L_{\omega}(\theta, \delta_0)\right) \, \Delta_{\omega}(\theta, g) \,,
%\\ \, & = \, & 
\end{eqnarray*}
which is indeed (\ref{Deltaellleq}), by virtue of Lemma \ref{lemmajjmp2006} and since $L_{\omega}(\theta, \delta_0) \, = \, (1-\omega) \, \|\delta_0 - \gamma(\theta)\|^2$.  \qed

A basic result for estimating a mean vector $\theta$  under quadratic loss,
for scale mixtures of normal distributions is the following.

\begin{lemma}
(Strawderman, 1974)
\label{lemma: straw74}
Let $X \sim f(\|x-\theta\|^2)$ have a scale mixture of normals distribution as in (\ref{representation-mixture}) with $d \geq 3$, and consider estimating $\theta$ with loss $\|\delta-\theta\|^2$.   Consider Baranchik-type estimators 
$\delta_{a,r(\cdot)}(X)$ as in (\ref{baranchik}) with conditions (\ref{conditions}).  Then, $\delta_{a,r(\cdot)}(X)$ dominates $X$ provided 
\begin{equation}
\nonumber
0 \,< a \, < \frac{2}{\mathbb{E}_0 (\frac{1}{\|X\|^2})}\,\,,
\end{equation}
and  provided $E_0[\|X\|^2]$  and $E_0[1/ ||X||^2]$ are finite.
\end{lemma}

The main result of this section can now be presented and established.

\begin{theorem}
\label{theoremdominanceforl}
Let $X \sim f(\|x-\theta\|^2)$ have a scale mixture of normals distribution as in (\ref{representation-mixture}) with $d \geq 3$, and consider estimating $\theta$ with loss $L_{\omega, \ell}$, as in (\ref{loss-section4}) with $\delta_0(X)=X$, where $\ell'$ satisfies Assumption \ref{assumptionell}.  Consider Baranchik-type estimator 
$\delta_{a(1-\omega),r(\cdot)}(X)$ as in (\ref{baranchik}) with conditions (\ref{conditions}) on $r(\cdot)$.   Then,  assuming that $f_0^*(x) \, = \, \ell'((1-\omega)\|x\|^2) \, f(\|x\|^2)/K_1$ is a density on $\mathbb{R}^d$,  $\delta_{a(1-\omega),r}(X)$  dominates $X$ provided 
\begin{equation}
\nonumber
0 \,< a \, < \frac{2}{\mathbb{E}_{0,\omega}^{*} (\frac{1}{\|X\|^2})}\,\,,
\end{equation}
and  provided $E_{0,\omega}^{*}[\|X\|^2]$  and $E_{0,\omega}^{*}[1/ ||X||^2]$ are finite,
where the expectation $\mathbb{E}_0^{*}$ is taken with respect to $f_0^*$.
\end{theorem}

{\bf Proof.}  With the given notation, observe that $\delta_{a(1-\omega),r}(X)= X + (1-\omega) \, g_{a,r}(X)$ with $g_{a,r}(X) = - \frac{a r(\|X\|^2)}{\|X\|^2}   \, X\,$.   Therefore, by Lemma \ref{lemmaDeltaforl} with $\delta_0(X)=X$ and $\gamma(\theta)=\theta$, we have for the difference in losses between $\delta_{a(1-\omega),r}(X)$ and $X$:
\begin{eqnarray*}
\Delta_{\omega,\ell}\,(\theta, g_{a,r}) \, & \leq &\, (1-\omega)^2 \,\ell' \left\lbrace 
(1-\omega) \, \|X - \theta \|^2 \right\rbrace  \, \Delta_0\,(\theta, g_{a,r})\,   
\\
 \, & = & \, K_1 \, (1-\omega)^2 \, \int_{\mathbb{R}^d } \, \left\lbrace \|\delta_{a,r}(x) - \theta \|^2 \, - \, \|x-\theta\|^2  \right\rbrace  f_0^*(\|x-\theta\|^2 \, dx \,.
\end{eqnarray*}
Finally, since both $\ell'$ and $f$ are completely monotone, so is the density $f_0^*$ (Lemma \ref{lemma-completelymonotone}, part a).  This implies that $f_0^*$ is a scale mixture of normals density and the result thus follows immediately from Lemma \ref{lemma: straw74}.
 \qed

\begin{remark}
For the unbalanced case $\omega=0$, one recovers Theorem 2.1 of Kubokawa et al. (2015).  For the original balanced loss function with $\ell(t)=t$, one may recover the result of Theorem \ref{theoremdominanceforl} directly by relying on Lemma \ref{lemma: straw74} and Lemma \ref{lemmajjmp2006}, as illustrated in Example \ref{examplesel} for the multivariate normal case.

Moreover, it is interesting to compare the balanced and unbalanced cut-off points $a_0(\omega)=\frac{2}{\mathbb{E}_{0,\omega}^{*} (\frac{1}{\|X\|^2})}$ and $a_0(0)=\frac{2}{\mathbb{E}_{0,0}^{*} (\frac{1}{\|X\|^2})}$.
For loss $L_{\omega,\ell}$ with $\ell(t)=t^{\beta}$, $0 < \beta <1 $, we have $a_0(\omega)=a_0(0)$ for all $\omega \in (0,1)$.  For choices $\ell$ such that $\frac{\ell'((1-\omega) t)}{ \ell'(t)}$ is non-decreasing in $t$,  we have a monotone likelihood ratio ordering for densities $f_0^*$ and $f$, with  the former being stochastically larger.   This implies the ordering $\mathbb{E}_{0,\omega}^{*} (\frac{1}{\|X\|^2}) \leq \mathbb{E}_{0,0}^{*} (\frac{1}{\|X\|^2})$, and therefore $a_0(\omega) \geq a_0(0)$.   Examples where such a condition holds include:  $\ell(t) = 1 - e^{-t/\alpha}$, $\alpha>0$, $\ell(t) = \log (1+t)$,  $\ell(t) = (1+t/\beta)^{\alpha}$ with $\beta>0, 0 < \alpha < 1$.
\end{remark}

It is also interesting to assess how the upper cut-off point on the multiple $a$ for the estimator $\delta_{a, r(\cdot)}$ to dominate $X$ varies in terms of the model $f$ and the choice of $\ell$ for the loss function.  In the former case, one can infer dominance results that are robust, holding for a given $f$ but also persisting for a class of departures from $f$.   This is quite plausible and simple to visualize as the cut-off point depends on $f$ only through the inverse moment 
$\mathbb{E}_{0,\omega}^{*} (\frac{1}{\|X\|^2})$.  For the latter case, one can infer dominance results that hold simultaneously for a subclass of losses (\ref{loss-section4}).  Here is such an illustration.

\begin{corollary}
\label{corollart-simultaneous-dominance}
Consider the context of Theorem \ref{theoremdominanceforl}, for a given loss $L_{\omega,\ell}$, and a Baranchik-type estimator $\delta_{a(1-\omega),r(\cdot)}(X)$ which satisfies the given requirements for dominance of $\delta_0(X)=X$.  Then, the dominance persists for the original balanced loss $L_{w,\ell_0}$ 
with $\ell_0(t)=t$ and $\delta_0(X)=X$.
\end{corollary}
{\bf Proof.}   It suffices to show that 
\begin{equation}
\label{inequality-last}
\mathbb{E}_{0}^{*} (\frac{1}{\|X\|^2}) \geq \mathbb{E}_{0} (\frac{1}{\|X\|^2}) \,,
\end{equation}
where the expectation $\mathbb{E}_{0}^{*}$ is taken with respect to the density $f_0^*$ given in Theorem \ref{theoremdominanceforl}, and where we define the expectation $\mathbb{E}_{0}$ as taken with respect to the density $f(\|x\|^2)$.  Now observe that the ratio of these densities, proportional to  $\ell'((1-\omega) \|x\|^2)$ is decreasing in $\|x\|^2$ by assumptions on $\ell$.  We thus have a monotone likelihood ratio in $\|x\|^2$ ordering between the densities and inequality (\ref{inequality-last}) follows since $1/\|x\|^2$ is decreasing in $\|x\|^2$.  \qed

\section{Concluding remarks}

For a multivariate normal distributed $ X \sim N_d(\theta, \sigma^2 I_d)$, and more generally for a scale mixture of normals model $X \sim f(\|x-\theta\|^2)$, we have provided shrinkage estimators of $\theta$ that improve on the benchmark estimator $\delta_0(X)$ as measured by the frequentist risk associated with balanced loss functions of the types (\ref{lossrho}) and (\ref{lossell}), and with completely monotone $\rho$ and $\ell$.   Much of the approach is unified with respect to the choices of $f$ and either $\rho$ or $\ell$ and the findings represent analytical extensions to the original balanced loss with either identity $\rho$ or $\ell$, unavailable up to now.

The findings in this paper do not cover cases with unknown scale such as observations generated from a $N_d(\theta, \sigma^2 I_d)$ with unknown $\sigma^2$, such as earlier results on the original balanced loss function (e.g., Chung et al. 1999; Zinodiny, 2014), but we expect that the techniques presented here should be useful to derive corresponding results for analogs of loss functions (\ref{lossrho}) and (\ref{lossell}).  Finally,  it would be most interesting and welcomed to obtain Bayesian estimators that either satisfy our conditions of dominance, or dominate the benchmark $\delta_0(X)=X$ under the set-up of Theorems \ref{dominancesection3} and \ref{theoremdominanceforl}.

\section*{Appendix}

\subsection*{Proof of Lemma \ref{lemmatechnical}, part (a)}

With $T/\alpha^2=\|Z+\frac{\theta}{\alpha}\|^2 \sim  \chi^2_d(\, \|\theta\|^2/\alpha^2)$ 
and the Poisson representation of the non-central $\chi^2$ distribution
(i.e.,  $T/\alpha^2 |K \sim \chi^2_{d+2K}\,,\, K \sim \hbox{Poisson}(\lambda= \|\theta\|^2/2\alpha^2)$ ), we have for $d \geq 3$ and $\theta \neq 0$
\begin{eqnarray}
\nonumber
\mathbb{E}_{\alpha}(\frac{1}{T}) \, & = & \frac{1}{\alpha^2} \, e^{-\lambda} \, \sum_{k \geq 0} \frac{\lambda^k}{k!} \; \mathbb{E} [ \frac{1}{\chi^2_{d+2k}}] \\
\nonumber
\, & = & \frac{1}{\alpha^2} \, e^{-\lambda} \, \sum_{k \geq 0} \frac{\lambda^{k}}{k!} \;
\frac{1}{d+2k-2} 
%\end{eqnarray}
%\end{document}
\\
\nonumber
\, & = & \frac{2}{\|\theta\|^2} \, e^{-\lambda} \, \sum_{k \geq 0} \frac{\lambda^{k+1}}{(k+1)!} \;
\frac{k+1}{d+2k-2} \\
\nonumber
\, & = & \frac{2}{\|\theta\|^2} \, e^{-\lambda} \, \sum_{k \geq 0} \frac{\lambda^{k}}{k!} \, U(k)  \\
\label{U}
\, & = & \,  \frac{2}{\|\theta\|^2} \, \mathbb{E}_{\lambda} \, U(K)\,,
\end{eqnarray}
with $U(K) \, = \, \frac{K}{d+2K-4} \, \mathbb{I}_{\mathbb{N}_+}(K)\,$.    Since $U(k)$ is increasing in $k \in \mathbb{N}$ for $d \geq 3$, since $\lambda$ is decreasing in $\alpha$,  and since the Poisson($\lambda$) distribution has increasing monotone likelihood ratio in $K$ with parameter $\lambda$, it follows from the above that $\mathbb{E}_{\alpha}(\frac{1}{T})$ is decreasing indeed in $\alpha$ for $d \geq 3$.  

\section*{Acknowledgements}
 
\'Eric Marchand's research is supported in part by the 
Natural Sciences and Engineering Research Council of Canada, and William 
Strawderman's research is partially supported by a grant from the Simons 
Foundation (\#418098).

\baselineskip=1.00\normalbaselineskip

\end{document}